\documentstyle[12pt]{article}

  \newcommand{\const}{\rm const}
  \newcommand{\Var}{\rm Var}

  \newcommand{\vraisup}{\rm vraisup}

  \newcommand{\argmin}{\rm argmin}
  
  \newcommand{\Ent}{\rm Ent}
  \newcommand{\Sub}{\rm Sub}

\begin{document}

   \begin{center}

 {\bf  Exponential confidence region based on the projection density estimate.} \par

\vspace{4mm}

 {\bf  Recursivity of these estimations. }\\

\vspace{5mm}

{\bf M.R.Formica,  E.Ostrovsky, and L.Sirota.}

  \end{center}

\vspace{4mm}

 Universit\`{a} degli Studi di Napoli Parthenope, via Generale Parisi 13, Palazzo Pacanowsky, 80132,
Napoli, Italy. \\

e-mail: mara.formica@uniparthenope.it \\

\vspace{4mm}

 \ Department of Mathematics and Statistics, Bar-Ilan University, \\
59200, Ramat Gan, Israel. \\

e-mail: eugostrovsky@list.ru\\

\vspace{4mm}

\ Department of Mathematics and Statistics, Bar-Ilan University, \\
59200, Ramat Gan, Israel. \\

e-mail: sirota3@bezeqint.net \\

\vspace{4mm}


\begin{center}

 \ {\bf Abstract}

\end{center}

\vspace{4mm}

  \hspace{3mm}  We investigate the famous Tchentzov's projection density statistical estimation
  in order to deduce the exponential decreasing  tail of distribution for the natural normalized deviation. \par
 \ We modify these estimations assuming the square integrability of estimated function,
 to make it recursive form, which is more convenient for applications, however  they have at
 the same speed of convergence as the for  the classical ones in the composite Hilbert space norm. \par

\vspace{4mm}

\begin{center}

 \ {\sc  Key words and phrases.}

 \end{center}

 \vspace{4mm}

 \hspace{3mm}  Projection statistical density estimation, measure and measurable space, Riesz - Fisher theorem, integrability
 and square integrability, probability space, random series and variables, sample, 
  orthogonal function and series, Hilbert space and  composite norm, orthonormal complete system, order,
 convergence, adaptivity, Fourier coefficients, sample, random variable (r.v.), density, consistency, ordinary and  adaptive
 projection estimation, slowly varying function,  recursivity, error, observations, optimal rule, triangle inequality.\par

 \vspace{5mm}

\section{Introduction.}

\vspace{5mm}

 \hspace{3mm} Let $ \  (X = \{x\}, \cal{B}, \mu)  \ $  be  measurable space equipped with sigma finite and separable measure $ \ \mu. \ $ Separability
 implies that the metric space $ \  \cal{B}  \ $ relative the distance function

 $$
 D(A_1, A_2) := \mu(A_1 \setminus A_2) + \mu(A_2 \setminus A_1)
 $$
is separable. \par

\vspace{4mm}

  \hspace{3mm}  Let  also $ \ \xi_i, \ i = 1,2,\ldots, \ \xi := \xi_1 \ $ be a  {\it sample}, i.e. a sequence  of independent  identical
distributed random variables (r.v.) defined on certain probability space $ \ (\Omega = \{\omega\}, \cal{V}, \ {\bf P} ) \  $
with expectation $ \ {\bf E} \ $ and variance $ \ \Var, \ $ having at the same non - negative integrable and in addition
 {\it square integrable}  unknown density function $ \ f = f(x), \ x \in X \ $ relative the measure $ \ \mu: \ $

$$
{\bf P} (\xi_i \in A) = \int_A f(x) \ \mu(dx), \ A \in \cal{B},
$$
such that

\begin{equation} \label{square integrable}
||f||^2(L_2(X,\mu)) =    ||f||^2 = ||f||^2_2 := \int_X f^2(x) \ \mu(dx) < \infty,
\end{equation}
and as ordinary

$$
f(x) \ge 0; \hspace{3mm} \int_X f(x) \ \mu(dx)  = ||f||_1 = 1.
$$

\vspace{5mm}

\hspace{3mm} {\bf  Our goal in this report is building of the  exponential decreasing confidence region in the Hilbert space norm
for the unknown density  as well as to offer and investigate the recursive projection statistical estimate of the density
function consistent in the mixed Hilbert space norm, having however at the same speed of convergence as the classical optimal projection
estimates. } \par

\vspace{5mm}

 \  Let us clarify briefly  the difference between ordinary and recursive estimations on the example of classical kernel density estimates.
The classical {\it kernel} Parzen - Rosenblatt estimate  $ \ f[PR]_n(x) \ $  of $ \ f(x), \ x \in R^d,\ d = 1,2,\ldots \ $  has a form

$$
f[PR]_n(x) := \frac{1}{n \ h^d(n)} \sum_{i=1}^n V\left( \ \frac{x - \xi_i}{h(n)} \ \right),
$$
where $ \ h(n) \to 0, \ n \to \infty. \ $  In contradiction, the so - called {\it recursive} Wolverton -  Wagner's estimate $ \ f[WW]_n(x) \ $ is
defined as follows

$$
f[WW]_n(x) := n^{-1} \sum_{i=1}^n \frac{1}{h(i)^d} V \left( \ \frac{x - \xi_i}{h(i)}   \ \right).
$$
 \ The second estimate has at the same asymptotical properties as the first one  when $ \ n \to \infty \ $   but it obey's the recursive form

$$
f[WW]_n(x) = \frac{n-1}{n} f[WW]_{n-1}(x) + \frac{1}{h^d(n)} V \left( \ \frac{x - \xi_n}{h(n)}  \ \right), \ n \ge 2,
$$
therefore it is more convenient for the practical using. \par

\begin{center}

 \vspace{5mm}

\section{Main result. Non - adaptive  projection estimations.}

\vspace{5mm}

\end{center}

 \hspace{3mm}  We will start from the classical {\it projection }  Tchentzov's density estimations, see  the classical works  \ \cite{Tchentzov 1}, \cite{Tchentzov 2}.
The more modern {\it adaptive} estimations will be discuss briefly  further. \par

\vspace{4mm}

\ Let  $ \ \{\phi_k\} = \{ \phi_k(x) \}, \ x \in X \ $ be certain {\it complete} orthonormal system  functions in the (separable)
 Hilbert space $ \ L_2 = L_2(X,\mu): \ $

\begin{equation} \label{orthonormality}
(\phi_k, \phi_l) = \int_X \phi_k(x) \ \phi_l(x) \ \mu(dx) = \delta_k^l, \ k,l = 1,2,\ldots,
\end{equation}
where $ \ \delta_k^l \ $ is  the Kroneker's symbol. \par

 \ We have been  assumed above that the density function is square integrable: $ \ f(\cdot) \ \in L_2; \ $ therefore in has a form
 (theorem of Riesz - Fisher)

\begin{equation} \label{expression}
f(x) = \sum_{k=1}^{\infty} c_k \ \phi_k(x),
\end{equation}
empirical Fourier coefficients, so that

\begin{equation} \label{rho N}
\lim_{N \to \infty} \rho(N) = 0, \hspace{3mm} \rho(N) = \rho[f](N)  \stackrel{def}{=} \sum_{k = N+1}^{\infty} c_k^2.
\end{equation}
 \ As long as

$$
c_k = \int_X \phi_k(x)\ f(x) \ \mu(dx) \ = \ {\bf E} \phi_k(\xi),
$$
the unbiased and  consistent as $ \ n \to \infty \ $ with probability one as well as in the mean square sense statistical estimate 
of the Fourier coefficients $ \ c_k(n) \ $ has a form

\begin{equation} \label{hat ck}
c_k(n) \stackrel{def}{=} n^{-1}  \sum_{i=1}^n  \phi_k(\xi_i).
\end{equation}

 \ The classical so - called {\it projection}  estimate of $ \ f(x) \ $  offered by N.N.Tchentzov \ \cite{Tchentzov 1}, \cite{Tchentzov 2}
has a form

\begin{equation} \label{Proj est}
f_{N,n}(x) \stackrel{def}{=} \sum_{k=1}^N c_k(n) \ \phi_k(x), \ N = 2,3,4, \ldots,
\end{equation}
 the number $ \ N \ $ is called {\it order} of this estimation. \par

\vspace{4mm}

 \ The quantity of this approximation may be measured by means of the mixed (anisotropic)  square Hilbert norm

\begin{equation} \label{quantity}
\Delta(N,n)  =  \Delta[f_{N,n},f]  \stackrel{def}{=} {\bf E} || f_{N,n}(\cdot) - f(\cdot)||^2_2.
\end{equation}

\vspace{4mm}

 \ {\bf Lemma 2.1.} Suppose that the orthonormal system $ \ \{\phi_k(\cdot)\} \ $ is uniformly bounded:

\begin{equation} \label{M condition}
M \stackrel{def}{=}    \sup_{k = 1,2,\ldots} \ \sup_{x \in X}  \phi_k^2(x) < \infty.
\end{equation}

 \ Then  $ \ \Delta(N,n) \le A(N,n),  \ $  where

\begin{equation} \label{Delta estim}
A(N,n) \stackrel{def}{=} \frac{M \ N}{n} + \rho(N).
\end{equation}

\vspace{3mm}

  {\bf Proof.} By virtue of orthonormality

$$
|| f_{N,n}(\cdot) - f(\cdot)||^2_2 =  \sum_{k=1}^N (c_k(n) - c_k)^2 + \sum_{k=N+1}^{\infty} c^2_k =
\sum_{k=1}^N (c_k(n) - c_k)^2 + \rho(N).
$$
 \ Further,

$$
{\bf E}  (c_k(n) - c_k)^2 = n^{-1} \Var  \{ \ \phi_k(\xi) \ \} \le n^{-1} {\bf E} [ \ \phi_k(\xi) \ ]^2 =
$$

$$
n^{-1}  \int_X \phi^2_k(x) \ f(x) \ \mu(dx) \le M \ n^{-1};
$$
following

\begin{equation} \label{Arb estim}
  {\bf E} || f_{N,n}(\cdot) - f(\cdot)||^2_2 \le M N/n + \rho(N)  = A(N,n),
\end{equation}

\vspace{3mm}

Q.E.D. \par

\vspace{4mm}

 \hspace{2mm}  Introduce the following important  {\it optimal}, more exactly, pseudo - optimal,
  variable (number of summands)

\begin{equation} \label{inportant N0}
N_0 = N_0(n)  = N_0[n,M,\rho] \stackrel{def}{=}  \argmin_{N = 2,3,\ldots, n} A(N,n),
\end{equation}
and correspondingly the optimal rate of convergence of  offered projection estimate is equal to 

\begin{equation}\label{minimal value}
A^*(n) \stackrel{def}{=} \min_N A(N,n) = A(N_0,n).
\end{equation}

\vspace{4mm}

 \ {\bf Example 2.1.} Let $ \ c_k^2 \asymp C_1 /k^{1 + \gamma}, \ C_j, \gamma = \const  \in (0,\infty). \ $ Then

 $$
 \rho(N) \asymp C_2 N^{-\gamma}, \hspace{3mm} N_0(n) \asymp C_3 N^{-\gamma}
 $$
and finally

$$
A^*(n) \asymp C_4 n^{-\gamma/(1 + \gamma)}, \ n \ge 1.
$$

\vspace{4mm}

 \ {\bf Example 2.2.} Let now  $ \ c_k^2 \asymp C_5 \exp( - C_6 \ k),  \ k = 1,2,\ldots; \ $ then

$$
 \rho(N) \asymp C_7 \exp(- C_8 \ N), \hspace{3mm} N_0(n) \asymp C_9 \ \ln n, \ n \ge 2,
$$

and finally in this case

$$
A^*(n) \asymp C_{10} \ln n/ n,  \ n \ge 2.
$$

\vspace{4mm}

 \  Ibragimov I.A. and Hasminskii R.Z. in \cite{Ibragimov} proved that this speed of convergence is optimal in the mixed
 Hilbert norm   $ \  L_2(\Omega \otimes X) \ $   sense   on the widely classes of estimated  functions   $ \ f. \ $  \par

\vspace{4mm}

 \ {\bf Lemma 2.2.}

 \begin{equation} \label{key restrictions}
 c^2_{N_0} \ge \frac{M}{n}; \hspace{4mm} c^2_{N_0 + 1} \le \frac{M}{n}, \hspace{3mm} n,N_0 \ge 2.
 \end{equation}

 \vspace{3mm}

 \ {\bf Proof.}  Both the relations in  (\ref{key restrictions}) follows immediately from the obvious inequalities

\begin{equation} \label{adjacent eq}
A(N_0 - 1,n) \le A(N_0,n) \le A(N_0+1,n), \hspace{3mm} N_0 = N_0(n).
\end{equation}

\vspace{4mm}

 \ {\it Briefly, by definition of notation: }

$$
 c^2_{N_0(n)} \stackrel{def}{\approx} \frac{M}{n}, \ n \to \infty.
$$

\vspace{4mm}

 \ Inversely, let the estimations

 \begin{equation} \label{L key restrictions}
 c^2_{L} \ge \frac{M}{n}; \hspace{4mm} c^2_{L + 1} \le \frac{M}{n}, \hspace{3mm} n \ge 2.
 \end{equation}
for some positive integer value $ \ L = L(n), \ n \ge 2 \ $ holds true.  Question: is the value $ \ L = L(n) \ $ actually
the optimal numbers $ \ N_0: \ L(n) = N_0(n)? \ $

\vspace{4mm}

\ {\bf Lemma 2.3.} Let the integer valued number $ \ L = L(n), \ n \ge 2 \ $ be a minimal value for which the relation (\ref{L key restrictions})
holds true. Let also the numerical non - negative sequence of square of coefficients $ \ c^2_k \ $  be monotonically non - increasing for sufficiently
greatest values $ \ k: \ $

\vspace{3mm}

\begin{equation} \label{momotonicity}
\exists k_0 = 2,3,\ldots \hspace{3mm}  \forall k \ge k_0  \Rightarrow  c^2_{k+1} \le c^2_k.
\end{equation}

\vspace{3mm}

 Then $ \ L(n) = N_0(n) = \argmin_N A(N,n). \ $ \par

\vspace{3mm}

 \ {\bf Proof.} Let us ground first of all the existence and non - triviality of the variable $ \ N_0 = N_0(n). \ $ Note that
 $ \ \lim_{N \to \infty} A(N,n) = \infty,   \  $ therefore  $ \ N_0(n) < \infty, \ n \ge 2. \ $ Further, denote by $ \ [m] = \Ent(m) \ $
the integer part the number $ \ m. \ $ We have

$$
A([\sqrt{n}], \ [\sqrt{n} ] ) \le   C \ \frac{1}{\sqrt{n}} + \rho[f]([\sqrt{n}]) \to 0, \ n \to \infty,
$$
as long as $ \ f \in L_2(X,\mu). \ $ Following

$$
\lim_{n \to \infty} A^*(n) =  \lim_{n \to \infty} A(N_0(n),n) = 0.
$$

\vspace{3mm}

 \ Ultimately, let the estimations (\ref{L key restrictions})  is valid for some positive integer value $ \ L = L(n), \ n \ge 2. \ $
Let also $ \  l = 1,2,\ldots \ $ be arbitrary positive integer number. We must ground the inequality

\begin{equation} \label{must ineq}
\rho(L) + \frac{ML}{n} \le \rho(L+l) + \frac{M(N+l)}{n},
\end{equation}
which is in turn quite equivalent to the next one

\begin{equation} \label{next one}
\sum_{j=1}^l c_{L + j}^2 \le l \frac{M}{n}, \ l = 1,2,\ldots;
\end{equation}
which follows immediately from the conditions of Lemma 2.3, because $ \ c^2_{L+j} \le c^2_L,  \ j = 1,2,\ldots . \ $ \par
 \ The case of the negative values $ \ l \ $ may be considered quite analogously. Namely,  one has

\begin{equation} \label{must left ineq}
\rho(L) + \frac{ML}{n} \le \rho(L-l) + \frac{M(L-l)}{n}, \ l = 1,2,\ldots, L - 1;
\end{equation}
which is in turn quite equivalent to the next inequality

\begin{equation} \label{next left one}
\sum_{j=1}^l c_{L - j}^2 \ge l \frac{M}{n}, \ l = 1,2,\ldots, L - 1.
\end{equation}

\vspace{3mm}

 \ {\bf Remark 2.1.} Note that the assertion of Lemma 2.3 remains true only under both  the relations (\ref{next one})
 and (\ref{next left one}). \par

\vspace{3mm}

\ {\bf Remark 2.2.} Let an integer non - random  value  $ \ K = K(n), n \ge 2, \ K(n) = 2,3, \ldots, n - 1 \ $  be such that
there exists an universal constant $ \ \Gamma  = \const \in [1, \infty) \ $ for which

$$
\forall n, \ N \in [2, n - 1] \ \Rightarrow \   A(K(n),n) \le \Gamma  \ A(N,n).
$$

 \ Then obviously

$$
{\bf E} || \ f(K(n),n)(\cdot)  - f(\cdot) \ ||_2^2 \le \Gamma \cdot A^*(n) \ -
$$
quasi - optimality. \par

\vspace{3mm}

\ {\bf Remark 2.3.} As long as in the practice the values $ \  \ A(N,n) \ $ are unknown, one can use instead of him
  its  statistical consistent estimations (adaptivity). Namely, define the variables

$$
\tau_n(N) \stackrel{def}{=} \sum_{k=N+1}^{2 N} c_k^2(n), \ n,N \ge 2, \ N \le n/2.
$$

 \ It is proved in  \cite{Bobrov}, see also \cite{Ostrovsky 0} - \cite{Ostrovsky 4} that as $ \ n \to \infty \ $
 with  probability one 

$$
\tau_n(N) \asymp A(N,n),
$$
 and following

$$
\min_{N \in [2, n/2]} \tau_n(N) \sim A^*(n), \hspace{3mm} \argmin_{N \in [2, n/2]} \tau_n(N) \sim N_0(n),
$$
also almost everywhere. \par

 \vspace{5mm}

 \section{Exponential confidence region based on the projection estimate.}

\vspace{5mm}

  \ Let us return to the arbitrary projection estimate   (\ref{Arb estim}). Define the  (non - negative) variable

\vspace{3mm}

\begin{equation} \label{confid estim}
\Delta = \Delta[f]_{N,n} \stackrel{def}{=}    \frac{n}{M \ N} \ \left[ \ || f_{N,n}(\cdot) - f(\cdot)||^2_2 - \rho(N) \ \right],
\end{equation}
and its tail function

\begin{equation} \label{Tail Delta}
T[\Delta] (t) \stackrel{def}{=} {\bf P}(\Delta > t), \ t \ge 1,
\end{equation}
which may be used by the building of the confidence region for the unknown  density function $ \ f(\cdot) \ $ in the
$ \ L_2(X) \ $ sense. Here $ \ 2 \le N \le n-2, \ n \ge 5; \ $ the value $ \ N \ $ may be optimal or not. \par

\vspace{5mm}

{\bf Theorem 3.1.} We conclude under formulated conditions

\vspace{3mm}

\begin{equation} \label{Tail estim}
T[\Delta] (t) \le \exp \left[ \ - t /  \left( \ e \ M \ \right)  \ \right], \ t \ge  e \ M.
\end{equation}

\vspace{5mm}

\ {\bf Proof.} We will use the known theory of Grand Lebesgue Spaces (GLS) of random variables
having the exponential decreasing tail of distributions, see e.g.
\cite{Ahmed Fiorenza Formica at all}, \cite{anatriellofiojmaa2015}, \cite{anatrielloformicaricmat2016},
\cite{Buldygin},  \cite{Capone1}, \cite{Capone2}, \cite{Kozachenko-Ostrovsky 1985}, \cite{Kozachenko-Ostrovsky-Sirota Jan2017},
\cite{Kozachenko-Ostrovsky-Sirota Oct2017}, \cite{Ostrovsky 0} etc. Concrete, let $ \ \psi = \psi(p), \ 1 \le p < b = \const \le \infty \ $
be measurable strictly positive:  $ \ \inf_p  \psi(p) > 0 \ $ function. By definition, the Banach rearrangement invariant  Grand Lebesgue Space (GLS)
$ \ G\psi \ $ consists on all the random variables (r.v.)  $ \ \{\xi\} \ $  having a finite norm

\begin{equation} \label{Gpsi norm}
||\xi|| G\psi \stackrel{def}{=} \sup_{ p \in (1,b)} \left\{ \ \frac{||\xi||L_{p,\Omega}}{\psi(p)} \ \right\} < \infty,
\end{equation}
where as ordinary

$$
||\xi||L_{p,\Omega} = ||\xi||_p := \left[ \  {\bf E} |\xi|^p  \ \right]^{1/p}, \ 1 \le p < \infty,
$$

$$
||\xi||_{\infty} = \vraisup_{\omega \in \Omega} |\xi(\omega)|.
$$

 \ The finiteness of some GLS norm $ \ ||\xi||G\psi \ $ for the r.v. $ \ \xi \ $ is closely related with its tail behavior

  $$
   T[\xi](t) \stackrel{def}{=} {\bf P}(|\xi| > t), \ t \ge e. \
  $$

Indeed, assume for definiteness that $ \ ||\xi||G\psi = 1; \ $ then

\begin{equation} \label{tail ineq}
T[\xi](t) \le \exp \left( \  - \sup_{p \in (1,b)} (p \ln t - p \ln \psi(p))\ \right), \ t \ge e;
\end{equation}
and inverse conclusion is true under simple  natural conditions: if (\ref{tail ineq}) there holds, then

$$
\xi \in G\psi, \ \hspace{3mm} \Leftrightarrow \ \exists  \ K = K(\psi) < \infty, \ ||\xi||G\psi \le K(\psi).
$$

 \vspace{3mm}

 \ {\bf Example.}  Set for certain $ \ m = \const \in (0,\infty) \ $

$$
\psi_m(p) := p^{1/m}, \ 1 \le p < \infty.
$$

 \ The GLS estimate

\begin{equation} \label{psi m}
||\xi||G\psi_m = \sup_{p \ge 1} \left[ \ \frac{||\xi||L(p,\Omega)}{\psi_m(p)}  \ \right] < \infty
\end{equation}
is completely equivalent to the following tail inequality

\begin{equation} \label{tail m}
\exists \gamma = \gamma(m) \in (0,\infty) \ \Rightarrow T[\xi](t) \le \exp \left( - \ \gamma(m) \ t^m \  \right),  \ t \ge 0.
\end{equation}

\vspace{3mm}

 \ The case $ \ m = 2 \ $ correspondent to the famous subgaussian random variables. Standard notation \ \cite{Kozachenko-Ostrovsky 1985}:

$$
||\xi||\Sub \stackrel{def}{=} ||\xi||G\psi_2 = \sup_{p \ge 1} \left[ \ \frac{||\xi||_p}{\sqrt{p}}  \ \right].
$$

 \ It is known that $ \ ||\xi||\Sub \le ||\xi||_{\infty}  \ $ and if the r.v. - s $ \ \{\eta_i\}, \ i = 1,2,\ldots,n  \  $ are subgaussian,
 independent, centered and identical distributed, then

\begin{equation} \label{sub estim}
\sup_n \ \left[ \ n^{-1/2} ||\sum_{i=1}^n \eta_i||\Sub \ \right] = ||\eta_1||\Sub.
\end{equation}

\vspace{4mm}

 \ Let us return to the theorem 3.1.  Denote

$$
\delta_k(n)  := n^{1/2} (c_k(n) - c_k) = n^{-1/2} \sum_{i=1}^n (\phi_k(\xi_i) - c_k);
$$

therefore

$$
n^{1/2} ||c_k(n) - c_k || \Sub \le \sqrt{M}.
$$

 \  We  entail  from here

$$
|| (c_k(n) - c_k)^2 ||G\psi_1 \le \frac{M}{n}.
$$

 \ One can apply the triangle inequality for the $ \ G\psi_1 \ $ norm

\begin{equation} \label{triangle}
|| \ \sum_{k=1}^N  (c_k(n) - c_k)^2 \ ||G\psi_1 \le \frac{M \ N}{n}.
\end{equation}

 \ Denote

$$
\Theta = \Theta[N,n] := \frac{n}{M \ N} \cdot  \sum_{k=1}^N  (c_k(n) - c_k)^2;
$$
then $ \  ||\Theta||G\psi_1 \le 1  \ $ and hence

$$
T[\Theta](t) \le \exp ( - t/e), \ t \ge e.
$$

 \ The proposition  (\ref{Tail estim}) follows immediately from the last estimate after simple calculations.\par

\vspace{5mm}

\begin{center}

 \section{Recursions. } \par

\vspace{5mm}

\end{center}

 \hspace{3mm} We intent in this section to offer a {\it recursive} modification of the projection density estimation
 having at the same up to  finite multiplicative constant speed of convergence as classical ones. \par

\vspace{3mm}

  \ Notice that the  number of summands $ \ N(n) \ $ is {\it integer positive number, } which is monotonically non -
  increasing relative the number $ \ n. \ $   Therefore  it is reasonable to suppose
that  it satisfied for all the sufficiently greatest values $ \ n, \ $ say for $ \ n \ge 10, \ $ the  following  recursion

\begin{equation} \label{n1}
N(n+1) = N(n) + 1
\end{equation}
or
 \begin{equation} \label{n0}
N(n+1) = N(n).
\end{equation}

 \ The first case (\ref{n1}) take place iff

\begin{equation} \label{grew case}
c^2_{N(n)} > \frac{M}{n},
\end{equation}
the second one  take place when

\begin{equation} \label{fix point}
c^2_{N(n)} \le \frac{M}{n}.
\end{equation}

 \vspace{3mm}

 \ So, let us introduce the following {\it recursively defined}  sequence of positive non - decreasing integer
 numbers $ \ Y = \{ Y(n) \} \ $ such that $ \  Y(1) = 1  \ $ and consequently

\begin{equation} \label{greating}
c^2_{Y(n)} > \frac{M}{n} \ \Rightarrow Y(n+1) := Y(n) + 1,
\end{equation}
and otherwise

\begin{equation} \label{stoping}
c^2_{Y(n)} \le \frac{M}{n} \ \Rightarrow Y(n+1) := Y(n).
\end{equation}

\vspace{5mm}

 \  {\bf Theorem 4.1.} We retain all the conditions and notations of the theorem 3.1.
   Define the following centered and normalized deviation variable

 \vspace{3mm}

\begin{equation} \label{confid  Y estim}
\Delta(Y) = \Delta[f]_{Y,N,n} \stackrel{def}{=}    \frac{n}{M \ Y(n)} \ \left[ \ || f_{Y(n),n}(\cdot) - f(\cdot)||^2_2 - \rho(Y(n)) \ \right],
\end{equation}

\vspace{3mm}

and its tail function

\vspace{3mm}

\begin{equation} \label{Tail Delta  Y}
T[\Delta(Y)] (t) \stackrel{def}{=} {\bf P}(\Delta(Y) > t), \ t \ge 1,
\end{equation}

\vspace{3mm}

 \ Here as above $ \ 2 \le Y(n) \le n-2, \ n \ge 5. \ $  \par

\vspace{3mm}

  Assume in addition that

 \vspace{3mm}

\begin{equation} \label{key restriction}
Q = Q[Y] \stackrel{def}{=} \sup_n \left\{ \ \frac{A(Y(n),n)}{A^*(n)} \ \right\} < \infty.
\end{equation}

\vspace{3mm}

 \ We conclude under formulated conditions by virtue of the proposition of remark 2.2

\vspace{3mm}

\begin{equation} \label{Tail estim}
T[\Delta(Y)] (t) \le \exp \left[ \ - t /  \left( \ e \ M \ Q \ \right)  \ \right], \ t \ge  e \ M \ Q.
\end{equation}

\vspace{5mm}

 \ {\bf Remark 4.1.} The condition   $ \ Q(Y) < \infty \ $ in  (\ref{key restriction}) is satisfied for instance when

$$
\exists C_1, C_2 \in (0,\infty), \ c^2_k \asymp k^{-1 - C_1} \ L(k), \ k \ge 2,
$$
where $ \ L = L(z), \ z \in [2,\infty) \ $ is non - negative continuous slowly varying at infinity function; as well as when

$$
\exists C_1, C_2 \in (0,\infty), \ c^2_k \asymp \exp \left( \ - C_1 k^{C_2} \ L(k)  \ \right),  \ k \ge 2.
$$

\vspace{5mm}

\section{Concluding remarks.}

\vspace{5mm}

 \hspace{3mm}  It is interest  in our opinion  to deduce the recursive  projection density (and regression, 
 spectral density) function estimates for the so - called {\it adaptive} estimates, having however  the optimal rate of 
 convergence, which does not dependent on the unknown, in general case,  class of  smoothness  for
 estimating function.  See for example   \cite{Bobrov}, \cite{Ostrovsky 1}, \cite{Ostrovsky 2}. \par

\vspace{6mm}

\vspace{0.5cm} \emph{Acknowledgement.} {\footnotesize The first
author has been partially supported by the Gruppo Nazionale per
l'Analisi Matematica, la Probabilit\`a e le loro Applicazioni
(GNAMPA) of the Istituto Nazionale di Alta Matematica (INdAM) and by
Universit\`a degli Studi di Napoli Parthenope through the project
\lq\lq sostegno alla Ricerca individuale\rq\rq .\par

\end{document}